\documentclass[11pt,a4paper]{amsart}
\usepackage[mathcal]{eucal}
\usepackage{color,graphicx}

\theoremstyle{plain}

\newtheorem{thmABC}{Theorem}

\newtheorem{corABC}[thmABC]{Corollary}

\newtheorem*{thm*}{Theorem}

\theoremstyle{remark}
\newtheorem*{acknowledgements}{Acknowledgements}

\newcommand{\N}{\mathbb{N}}
\newcommand{\Z}{\mathbb{Z}}
\newcommand{\Q}{\mathbb{Q}}
\newcommand{\C}{\mathbb{C}}
\newcommand{\R}{\mathbb{R}}

\newcommand{\Shadows}{\mathfrak{S}\mathfrak{h}}
\newcommand{\lri}{\mathfrak{o}}
\newcommand{\Lri}{\mathfrak{O}}
\newcommand{\gri}{\ensuremath{{\scriptstyle \mathcal{O}}}}
\newcommand{\smallgri}{\ensuremath{{\scriptscriptstyle \mathcal{O}}}}

\renewcommand{\epsilon}{\varepsilon}

\DeclareMathOperator{\rk}{rk}

\DeclareMathOperator{\SL}{SL}
\DeclareMathOperator{\GL}{GL}

\DeclareMathOperator{\SU}{SU}
\DeclareMathOperator{\Tr}{Tr}
\DeclareMathOperator{\real}{Re}
\DeclareMathOperator{\sh}{sh}

\DeclareMathOperator{\Stab}{Stab}
\newcommand{\nir}[1]{#1}

\begin{document}
\title[Representation zeta functions of groups]{On representation zeta functions
  of groups \\ and a conjecture of Larsen and Lubotzky}

\author{Nir Avni}\thanks{\nir{Avni was supported by NSF grant
  DMS-0901638.}} \address{Department of Mathematics, Harvard
  University, One Oxford Street, Cambridge MA 02138, USA}
  \email{avni.nir@gmail.com}

\author{Benjamin Klopsch} \address{Department of Mathematics, Royal
  Holloway, University of London, Egham TW20 0EX, United Kingdom}
\email{Benjamin.Klopsch@rhul.ac.uk}

\author{Uri Onn} \address{Department of Mathematics, Ben Gurion
  University of the Negev, Beer-Sheva 84105 Israel}
\email{urionn@math.bgu.ac.il}

\author{Christopher Voll} \address{School of Mathematics, University
  of Southampton, University Road, Southampton SO17 1BJ, United
  Kingdom} \email{C.Voll.98@cantab.net}

\begin{abstract}

We study zeta functions enumerating finite-dimensional irreducible
complex linear representations of compact $p$-adic analytic and of
arithmetic groups. Using methods from $p$-adic integration, we
show that the zeta functions associated to certain $p$-adic
analytic pro-$p$ groups satisfy functional equations. We prove a
conjecture of Larsen and Lubotzky regarding the abscissa of
convergence of arithmetic groups of type $A_2$ defined over number
fields, assuming a conjecture of Serre on lattices in semisimple
groups of rank greater than $1$.

\end{abstract}

%\thanks{\hfill \textbf{Date of draft version: \today}}


% {analysis on p-adic lie groups, discrete subgroups of lie groups,
%   other zeta and Dirichlet functions, ordinary representations and
%   characters, linear algebraic groups over local fields} {20E18,
%   20C15, 20G25, 20E50}

\maketitle

%%%%%

\section{Introduction and statement of results}
Let $G$ be a group and, for $n \in \N$, denote by $r_n(G)$ the number
of equivalence classes of $n$-dimensional irreducible complex
representations of $G$; if $G$ is a topological or an algebraic group,
it is tacitly understood that representations are continuous or
rational, respectively. We assume henceforth that $G$ is
(representation) \emph{rigid}, i.e.\ that $r_n(G)$ is finite for
all~$n\in\N$.  In the subject of representation growth one
investigates the arithmetic properties of the sequence $r_n(G)$ and
its asymptotic behaviour as $n$ tends to infinity.  Recent key
advances in this area were made by Larsen and Lubotzky in
\cite{LaLu08}.

The group $G$ is said to have \emph{polynomial representation
growth} (PRG) if the sequence $R_N(G):=\sum_{n=1}^Nr_n(G)$ is
bounded by a polynomial.  An important tool to study the
representation growth of a PRG group $G$ is its
\emph{representation zeta function}, viz.\ the Dirichlet series
$$
\zeta_{G}(s):=\sum_{n=1}^\infty r_n(G)n^{-s},
$$ where $s$ is a complex variable.  It is well-known that the
\emph{abscissa of convergence} $\alpha(G)$ of the series $\zeta_G(s)$,
i.e.\ the infimum of all $\alpha\in\R$ such that $\zeta_G(s)$
converges on the complex right half-plane $\{s \in \C \mid \real (s)>
\alpha \}$, gives the precise degree of polynomial growth:
$R_N(G)=O(1+N^{\alpha(G)+\epsilon})$ for every $\epsilon \in \R_{>0}$.

In \cite{AvKlOnVo09I} we introduce new methods from the theory of
$\mathfrak{p}$-adic integration to study representation zeta functions
associated to compact $p$-adic analytic groups and arithmetic
groups. In \cite{AvKlOnVo09II} we compute explicit formulae for the
representation zeta functions of the groups $\SL_3(\lri)$, where
$\lri$ is a compact discrete valuation ring of characteristic $0$, in
the case that the residue field characteristic is large compared to
the ramification index of~$\lri$.  We give a summary of the main
results of our forthcoming papers in the current section, followed by
a brief description of the methodology in
Section~\ref{sec:methodology}.

\smallskip

A finitely generated profinite group $G$ is rigid if and only if it is
FAb, i.e.\ if every open subgroup of $G$ has finite
abelianisation. In~\cite{Ja06}, Jaikin-Zapirain proved rationality
results for the representation zeta functions of FAb compact $p$-adic
analytic groups using tools from model theory. In particular, the
representation zeta function of a FAb $p$-adic analytic pro-$p$ group
is a rational function in $p^{-s}$, for $p>2$.  Key examples of FAb
compact $p$-adic analytic groups are the special linear groups
$\SL_n(\lri)$ and their principal congruence
subgroups~$\SL_n^m(\lri)$, where $\lri$ is a compact discrete
valuation ring of characteristic $0$ and residue field characteristic
$p$.  For fixed $n$, and varying $m$ and $\lri$, the latter also yield
important examples of families of pro-$p$ groups which arise from a
global Lie lattice, in this case $\mathfrak{sl}_n(\Z)$.

To be more precise, let $\gri$ be the ring of integers of a number
field~$k$, and let $\Lambda$ be an $\gri$-Lie lattice such that
$k\otimes_\smallgri \Lambda$ is a perfect $k$-Lie algebra of
dimension~$d$. Let $\lri=\gri_v$ be the ring of integers of the
completion $k_v$ of $k$ at a non-archimedean place~$v$, lying above a
rational prime~$p$.  Given a finite extension $\Lri$ of $\lri$, we
write $\mathfrak{P}$ for the maximal ideal of $\Lri$,
$e(\Lri\vert\lri)$ for the ramification index and $f(\Lri\vert\lri)$
for the residue class field extension degree. Let $\mathfrak{g}(\Lri)
:= \Lri\otimes_\smallgri \Lambda$. For all sufficiently large $m$, the
Lie lattice $\mathfrak{g}^m(\Lri) :=
\mathfrak{P}^m\,\mathfrak{g}(\Lri)$ corresponds, by $p$-adic Lie
theory, to a FAb, potent, saturable pro-$p$ group $\mathsf{G}^m(\Lri)
:= \exp \left( \mathfrak{g}^m(\Lri) \right)$.  We call such $m$
\emph{permissible} for the Lie lattice $\mathfrak{g}(\Lri)$.  \nir{For
example}, for unramified extensions $\Lri$ of $\Z_p$ and $p$ odd,
every $m \in \N$ is permissible.  In \cite{AvKlOnVo09I}, we prove

\begin{thmABC}%[{\cite[Theorem A]{AvKlOnVo09I}}]
  \label{thmABC:funeq}
  In the above setup, there exist a finite set $S$ of places of $k$,
  a natural number $r$ and a rational function
  $W(X_1,\dots,X_r,Y)\in\Q(X_1,\dots,X_r,Y)$ such that, for every
  non-archimedean place $v$ of $k$ with $v\not\in S$, the following is
  true.

  There exist algebraic integers
  $\lambda_1=\lambda_1(v),\dots,\lambda_r=\lambda_r(v)$ such that for
  all finite extensions $\Lri$ of $\lri=\gri_v$ and for all $m\in\N$
  which are permissible for $\mathfrak{g}(\Lri)$ one has
  \begin{equation}\label{equ:rational}
    \zeta_{\mathsf{G}^m(\Lri)}(s)=q_v^{fdm}\,W(\lambda_1^f,\dots,\lambda_r^f,q_v^{-fs}),
  \end{equation}
  where $q_v$ denotes the residue field cardinality of $\lri$, $f =
  f(\Lri\vert\lri)$ and $d=\mathrm{rank}_\Lri
  \left(\mathfrak{g}(\Lri)\right)=\dim_k(k\otimes_\gri \Lambda)$.
  Furthermore, the functional equation
   \begin{equation}\label{equ:funeq}
   \zeta_{\mathsf{G}^m(\Lri)}(s)|_{\substack{q_v\rightarrow
      q_v^{-1}\\\lambda_i\rightarrow
      \lambda_i^{-1}}}=q_v^{fd(1-2m)}\zeta_{\mathsf{G}^m(\Lri)}(s)
  \end{equation}
holds.
\end{thmABC}

Our proof of Theorem~\ref{thmABC:funeq} implies in particular that the
real parts of the poles of the zeta functions
$\zeta_{\mathsf{G}^m(\Lri)}(s)$ are rational numbers. More precisely,
we prove the following.
%\footnote{Nir
%proposed an argument that the monotonicity of the abscissa of
%convergence also holds when $v \in S$, but this seems to rely on
%extending some algebraic geometry results to schemes over rings?  Do
%we have a detailed proof?}

\begin{thmABC}% [{\cite[Theorem B]{AvKlOnVo09I}}]
  \label{thmABC:poles}
  In the above setup, there exists a finite set $P \subset \Q$ such
  that for all non-archimedean places $v$ of $k$, all finite
  extensions $\Lri$ of $\lri = \gri_v$, and all permissible $m$ for
  $\mathfrak{g}(\Lri)$ one has
  $$
  \{\real(s)\mid s \text{ a pole of }\zeta_{\mathsf{G}^m(\Lri)}(s) \}
  \subseteq P.
  $$

  Furthermore, if $v \not \in S$, where $S$ is a finite set of places
  arising from Theorem~\ref{thmABC:funeq}, and if $\gri_v \subseteq
  \Lri_1\subseteq \Lri_2$, then for every $m\in\N$ which is permissible
  for $\mathfrak{g}(\Lri_{1})$ and $\mathfrak{g}(\Lri_{2})$,
  \begin{equation}\label{equ:comparison_abscissae}
    \alpha(\mathsf{G}^m(\Lri_{1})) \leq \alpha(\mathsf{G}^m(\Lri_{2})).
  \end{equation}
\end{thmABC}

Notice that, if the groups $\mathsf{G}^m(\Lri)$ are principal
congruence subgroups of a FAb compact $p$-adic analytic group
$\mathsf{G}(\Lri)$ consisting of the $\Lri$-points of an algebraic
group $\mathbf{G}$, such as $\mathbf{G} = \SL_n$, then
\eqref{equ:comparison_abscissae} implies the monotonicity of the
abscissae of convergence $\alpha(\mathsf{G}(\Lri))$ under ring
extensions. This follows from the fact that these abscissae are
commensurability invariants. The set $P$ of candidate poles is
obtained by means of a resolution of singularities which leads to the
generic formula \eqref{equ:rational}.  Theorems~\ref{thmABC:funeq}
and~\ref{thmABC:poles} are illustrated by the explicit formulae given
in Theorem~\ref{thmABC:SL3} below.

The arithmetic groups we are interested in are arithmetic lattices
in semisimple algebraic groups defined over number fields. More
precisely, let $\mathbf{G}$ be a connected, simply connected
semisimple algebraic group, defined over a number field $k$,
together with a fixed $k$-embedding into some $\GL_n$.  Let
$\gri_S$ denote the ring of $S$-integers in~$k$, for a finite set
$S$ of places of $k$ including all the archimedean ones.  We
consider groups $\Gamma$ which are commensurable to
$\mathbf{G}(\gri_S)=\mathbf{G}(k) \cap \GL_n(\gri_S)$.
In~\cite{LuMa04}, Lubotzky and Martin showed that such a group
$\Gamma$ has PRG if and only if $\Gamma$ has the Congruence
Subgroup Property (CSP).  Suppose that $\Gamma$ has these
properties.  Then, according to a result of Larsen and Lubotzky
\cite[Proposition~1.3]{LaLu08}, the representation zeta function
of $\Gamma$ admits an Euler product decomposition. Indeed, if
$\Gamma = \mathbf{G}(\gri_S)$ and if the congruence kernel of
$\Gamma$ is trivial, this decomposition is particularly easy to
state: it takes the form
\begin{equation}\label{equ:euler}
  \zeta_\Gamma(s) = \zeta_{\mathbf{G}(\C)}(s)^{[k:\Q]} \cdot \prod_{v\not\in
    S}\zeta_{\mathbf{G}(\smallgri_v)}(s).
\end{equation}
Here each archimedean factor $\zeta_{\mathbf{G}(\C)}(s)$ enumerates
rational representations of the group $\mathbf{G}(\C)$; their
contribution to the Euler product reflects Margulis
super-rigidity. The groups $\mathbf{G}(\gri_v)$ are FAb compact
$p$-adic analytic groups whose principal congruence subgroups fit into
the framework of Theorems~\ref{thmABC:funeq} and~\ref{thmABC:poles};
the product of the zeta functions of these local groups captures the
finite image representations of $\Gamma$.

Several of the key results of \cite{LaLu08} concern the abscissae
of convergence of the `local' representation zeta functions
occurring as Euler factors on the right hand side
of~\eqref{equ:euler}. With regards to abscissae of convergence of
the `global' representation zeta functions Avni proved that, for
an arithmetic group $\Gamma$ with the CSP, the abscissa of
convergence of $\zeta_\Gamma(s)$ is always a rational number;
see~\cite{Av08}. In \cite[Conjecture 1.5]{LaLu08}, Larsen and
Lubotzky conjectured that, for any two irreducible lattices
$\Gamma_1$ and $\Gamma_2$ in a higher-rank semisimple group $H$,
one has $\alpha(\Gamma_1)=\alpha(\Gamma_2)$, i.e.\ that the
abscissa of convergence only depends on the ambient group. This
can be regarded as a refinement of Serre's conjecture on the
Congruence Subgroup Property. In \cite[Theorem~10.1]{LaLu08},
Larsen and Lubotzky prove their conjecture in the case that $H$ is
a product of simple groups of type $A_1$, assuming Serre's
conjecture. In~\cite{AvKlOnVo09I}, we prove

\begin{thmABC}%[{\cite[Theorem C]{AvKlOnVo09I}}]
  \label{thmABC:lalu}
  Let $\Gamma$ be an arithmetic \nir{lattice} of a connected, simply
  connected simple algebraic group of type $A_2$ defined over a number
  field.  If $\Gamma$ has the CSP, then $\alpha(\Gamma)=1$.
\end{thmABC}

\begin{corABC} Assuming Serre's conjecture, Larsen and Lubotzky's
  conjecture holds for groups of the form $H = \prod_{i=1}^r \mathbf{G}_i(K_i)$,
  where each $K_i$ is a local field of characteristic $0$ and each
  $\mathbf{G}_i$ is an absolutely almost simple $K_i$-group of type
%$A_1$ or
  $A_2$ such that $\sum_{i=1}^r \rk_{K_i}(\mathbf{G}_i) \geq 2$ and
  none of the $\mathbf{G}_i(K_i)$ is compact.
\end{corABC}

Key to our proof of Theorem~\ref{thmABC:lalu} is the following local
result, which we formulate in accordance with the notation introduced
before Theorem \ref{thmABC:funeq}.

\begin{thmABC}%[{\cite[Theorem D]{AvKlOnVo09I}}]
  \label{thmABC:SL3}
  Let $\lri$ be a compact discrete valuation ring of characteristic~$0$,
  with residue field of cardinality~$q$.  Let $\mathfrak{g}(\lri)$ be
  one of the following two $\lri$-Lie lattices of type $A_2$:
  \begin{enumerate}
  \item [(a)] $\mathfrak{sl}_3(\lri) = \{ \mathbf{x} \in
    \mathfrak{gl}_3(\lri) \mid \Tr(\mathbf{x}) = 0 \}$;
  \item [(b)] $\mathfrak{su}_3(\Lri,\lri) = \{ \mathbf{x} \in
    \mathfrak{sl}_3(\Lri) \mid \mathbf{x}^\sigma = -
    \mathbf{x}^{\mathrm{t}} \}$, where $\Lri \vert \lri$ is an
    unramified quadratic extension with nontrivial automorphism
    $\sigma$.
  \end{enumerate}

  For $m \in \N$, let $\mathsf{G}^m(\lri)$ be the $m$-th principal
  congruence subgroup of the corresponding group $\SL_3(\lri)$ or
  $\SU_3(\Lri,\lri)$. Assume that the residue field characteristic of
  $\lri$ is not equal to $3$. Then, for all $m \in \N$ which are
  permissible for $\mathfrak{g}(\lri)$, one has
  \begin{equation*}
    \zeta_{\mathsf{G}^m(\lri)}(s) = q^{8m} \frac{1 + u(q) q^{-3-2s}
      + u(q^{-1}) q^{-2-3s} + q^{-5-5s}}{(1 - q^{1-2s})(1 - q^{2-3s})},
  \end{equation*}
  where
  \begin{equation*} u(X) =
    \begin{cases}
      \phantom{-}X^3 + X^2 - X - 1 - X^{-1} & \text{ if }
      \mathfrak{g}(\lri)=
      \mathfrak{sl}_3(\lri),\\
      -X^3 + X^2 - X + 1 - X^{-1} & \text{ if }\mathfrak{g}(\lri)=
      \mathfrak{su}_3(\Lri,\lri).
    \end{cases}
  \end{equation*}
\end{thmABC}

The close resemblance between the representation zeta functions of
the special linear and the special unitary groups is noteworthy
and reminiscent of the Ennola duality for the characters of the
corresponding finite groups of Lie type.  We also give an explicit
formula for $\zeta_{\SL_3^m(\lri)}(s)$ in the exceptional case
where $\lri$ is unramified and has residue field
characteristic~$3$. Note that Theorem~\ref{thmABC:SL3} implies
that the abscissae of convergence $\alpha(\SL_3(\lri))$ and
$\alpha(\SU_3(\Lri,\lri))$ are each equal to $2/3$, as the
abscissa of convergence is a commensurability invariant.

The explicit formula for $\zeta_{\SL_3(\lri)}(s)$, which we present in
\cite{AvKlOnVo09II}, is deduced by means of the Kirillov orbit method,
a description of the similarity classes in finite quotients of
$\mathfrak{gl}_3(\lri)$, and Clifford theory.

\begin{thmABC}%[{\cite[Theorem 1.3]{AvKlOnVo09II}}]
  \label{zeta.SL3}
  There exist finitely many polynomials $f_{\tau, i},g_{\tau, i}\in
  \Q[x]$, indexed by $(\tau,i) \in \{1,-1\} \times I$, such that for
  every compact discrete valuation ring $\lri$ of characteristic $0$,
  with residue field characteristic $p > 3e(\lri \vert \Z_p)$, one has
  \[
  \zeta_{\SL_3(\lri)}(s)= \frac{\sum_{i \in I}
    f_{\tau,i}(q)(g_{\tau,i}(q))^{-s}}{(1 - q^{1-2s})(1 - q^{2-3s})},
  \]
  where $q$ denotes the size of the residue field of $\lri$ and $q
  \equiv \tau \pmod{3}$.
\end{thmABC}

An explicit set of such polynomials $f_{\tau, i},g_{\tau, i}$ is
computed in~\cite{AvKlOnVo09II}.  This result should be seen against
the background of ~\cite[Theorem~1.1]{Ja06}, which establishes the
rationality of representation zeta functions of FAb compact $p$-adic
analytic groups.  For groups of the form $\SL_3(\lri)$,
Theorem~\ref{zeta.SL3} specifies that these rational functions vary
`uniformly' as a function of the residue field
cardinality~$q$. Theorem~\ref{zeta.SL3} enables us to analyse the
global representation zeta functions of the arithmetic groups
$\SL_3(\gri_S)$.

\begin{thmABC}\label{thmABC:globalSL3}
  Let $\gri_S$ be the ring of $S$-integers of a number field $k$,
  where $S$ is a finite set of places of $k$ including all the
  archimedean ones.  Then there exists $\epsilon > 0$ such that the
  representation zeta function of $\SL_3(\gri_S)$ admits a meromorphic
  continuation to the half-plane $\{s \in \C \mid \real(s) >
  1-\epsilon\}$.  The continued function is holomorphic on the line
  $\{ s \in \C \mid \real(s)=1 \}$ except for a double pole at $s=1$.
  There is a constant $c \in \mathbb{R}_{>0}$ such that
  $$
  \sum_{n=1}^N r_n(\SL_3(\gri_S)) \sim c \cdot N \log N,
  $$
  where $f(N)\sim g(N)$ means $\lim_{N \rightarrow \infty}
  f(N)/g(N) = 1$.
\end{thmABC}

In~\cite{AvKlOnVo09I}, we also give simple geometric estimates for the
abscissae of convergence of representation zeta functions of compact
$p$-adic analytic groups and we compute representation zeta functions
associated to norm-$1$ groups in non-split quaternion algebras.

\section{Methodology}\label{sec:methodology}
\subsection{Kirillov orbit method and $\mathfrak{p}$-adic integration}

The core technique in \cite{AvKlOnVo09I} is a $\mathfrak{p}$-adic
formalism for the representation zeta functions of potent, saturable
pro-$p$ groups.  This approach has two key ingredients. Firstly, the
\emph{Kirillov orbit method} for potent, saturable pro-$p$ groups
provides a way to construct the characters of these groups in terms of
co-adjoint orbits; cf., e.g.,~\cite{Go08}. This `linearisation' --
pioneered in~\cite{Ho77b,Ja06} -- allows us to transform the original
problem of enumerating representations by their dimension into the
problem of counting co-adjoint orbits by their size. The second main
idea of our approach is to tackle the latter problem with the help of
suitable \emph{$\mathfrak{p}$-adic integrals} which are closely
related to Igusa local zeta functions (cf.~\cite{De91, Ig00, VeZu08})
and conceptually simpler than the definable integrals utilised
in~\cite{Ja06}.

For pro-$p$ groups of the form $\mathsf{G}^m(\Lri)$, which arise from
a global $\gri$-Lie lattice $\Lambda$ as in the setup of
Theorems~\ref{thmABC:funeq} and~\ref{thmABC:poles}, we describe the
representation zeta functions in terms of $\mathfrak{p}$-adic integrals of the
shape
\begin{equation} \label{equ:integral_intro} \mathcal{Z}_{\Lri}(r,t) =
  \int_{(x,\mathbf{y}) \in V(\Lri)} \lvert x \rvert_\mathfrak{P}^t
  \prod_{1 \leq j \leq \lfloor d/2 \rfloor} \frac{\lVert
    F_j(\mathbf{y}) \cup F_{j-1}(\mathbf{y})x^2
    \rVert_\mathfrak{P}^r}{\lVert F_{j-1}(\mathbf{y})
    \rVert_\mathfrak{P}^r} \, d\mu(x,\mathbf{y}).
\end{equation}
Here $\mathfrak{P}$ denotes the maximal ideal of the discrete
valuation ring $\Lri$, the domain of integration
$V(\Lri)\subset\Lri^{d+1}$ is a union of cosets modulo $\mathfrak{P}$,
the additive Haar measure $\mu$ is normalised so that
$\mu(\Lri^{d+1})=1$, the $F_j(\mathbf{y})$ are families of polynomials
over the global ring~$\gri$, which may be defined in terms of the
structure constants of the lattice $\Lambda$ with respect to a given
$\gri$-basis, and we write $\Vert \cdot\Vert_{\mathfrak{P}}$ for
the $\mathfrak{P}$-adic maximum norm.

%\nir{and we denote the maximum norm of an element in the set $S$ by $\Vert S\Vert_{\mathfrak{P}}$.}

The link between the Kirillov orbit formalism and the $\mathfrak{p}$-adic
integrals \eqref{equ:integral_intro} is given by the fact that the
problem of enumerating finite co-adjoint orbits of given size may be
reformulated as the problem of enumerating elementary divisors of
matrices of linear forms.  The integrals~\eqref{equ:integral_intro}
are multivariate analogues of Igusa local zeta functions associated to
polynomial mappings, a well-studied class of local zeta integrals;
see~\cite{VeZu08}. They are also akin to the $\mathfrak{p}$-adic integrals
studied in~\cite{Vo}, and our proofs of Theorems~\ref{thmABC:funeq}
and~\ref{thmABC:poles} rely on an adaptation of the methods and
formulae provided there.

A key point in the proof of Theorem~\ref{thmABC:lalu} is the fact
that, for the relevant arithmetic groups $\Gamma$, almost all of the
non-archimedean Euler factors in~\eqref{equ:euler} are of the form
$\zeta_{\SL_3(\lri)}(s)$ or $\zeta_{\SU_3(\Lri,\lri)}(s)$, and that
for the exceptional factors the abscissa of convergence is
sufficiently small, in fact equal to $2/3$. In \cite{AvKlOnVo09I} we
use the exact formulae provided by Theorem~\ref{thmABC:SL3}, together
with Clifford theory and suitable approximative Dirichlet series, to
prove that the global abscissa of convergence $\alpha(\Gamma)$ is
always equal to~$1$.

The proof of Theorem~\ref{thmABC:SL3} relies on a concrete
interpretation of the $\mathfrak{p}$-adic
integrals~\eqref{equ:integral_intro}. More generally, in the case of
`semisimple' compact $p$-adic analytic groups, we give a description
of these integrals in terms of a filtration of the irregular locus of
the associated Lie algebra. The terms of this filtration are
projective algebraic varieties defined by centraliser dimension,
refining the notion of irregularity. In the case of $\mathfrak{sl}_3$
and $\mathfrak{su}_3$, this filtration is simple enough to allow
explicit computations.

\subsection{Similarity classes, shadows and Clifford
  theory}\label{sec:p2}

The explicit computation of the zeta function of $\SL_3(\lri)$ in
Theorem~\ref{zeta.SL3} is based on the Kirillov orbit method and a
quantitative analysis of the $\GL_3(\lri)$-adjoint orbits, or
\emph{similarity classes}, in the finite quotients
$\mathfrak{gl}_3(\lri/\mathfrak{p}^\ell)$, $\ell\in\N$. The
pre-images of similarity classes in
$\mathfrak{gl}_3(\lri/\mathfrak{p}^\ell)$ under the natural
reduction map
$\pi_\ell:\mathfrak{gl}_3(\lri/\mathfrak{p}^{\ell+1})\to\mathfrak{gl}_3(\lri/\mathfrak{p}^\ell)$
are unions of similarity classes. In~\cite{AvKlOnVo09II} we give
an explicit, recursive description of these classes which allow us
to compute their numbers and cardinalities.

%Understanding the decomposition of
%these pre-images gives an inductive process through which the numbers
%and cardinalities of the orbits can be computed.

Given $\ell\in\N$ and a similarity class $\mathcal{C} \in
\GL_3(\lri)\backslash\mathfrak{gl}_3(\lri/\mathfrak{p}^\ell)$, we
define the shadow $\sh(\mathcal{C})$ of $\mathcal{C}$ to be the
conjugacy class of the reduction modulo $\mathfrak{p}$ of
$\Stab_{\GL_3(\lri)}(A)$, where $A$ is any element of
$\mathcal{C}$. This definition is independent of the choice of $A$. We
write $\Shadows:=\{\sh(\mathcal{C})\mid\mathcal{C} \in
\GL_3(\lri)\backslash\mathfrak{gl}_3(\lri/\mathfrak{p}^\ell),
\ell\in\N\}$ for the set of shadows. It turns out that
$\vert\Shadows\vert=10$, and that the shadow of a similarity class
$\mathcal{C}$ determines the numbers and sizes of the classes which
make up $\pi_\ell^{-1}(\mathcal{C})$. More precisely, there are
polynomials
$a_{\sigma_1,\sigma_2}(x),b_{\sigma_1,\sigma_2}(x)\in\Q[x]$, indexed
by $(\sigma_1,\sigma_2)\in\Shadows^2$, such that, for all $\ell\in\N$,
any similarity class $\mathcal{C} \in
\GL_3(\lri)\backslash\mathfrak{gl}_3(\lri/\mathfrak{p}^\ell)$ and
shadow $\sigma\in\Shadows$, the following hold: %if
%$\mathcal{D}$ is the pre-image of $\mathcal{C}$ under the reduction
%map
%$\mathfrak{gl}_3(\lri/\mathfrak{p}^{\ell+1})\to\mathfrak{gl}_3(\lri/\mathfrak{p}^\ell)$,
%then

\begin{enumerate}
\item [(a)] The number of similarity classes in
$\pi_\ell^{-1}(\mathcal{C})$ with shadow $\sigma$ is equal to
$a_{\sh(\mathcal{C}),\sigma}(q)$.
\item [(b)] All similarity classes in $\pi_\ell^{-1}(\mathcal{C})$
with shadow $\sigma$ have equal cardinality $|\mathcal{C}|\cdot
b_{\sh(\mathcal{C}),\sigma}(q)$.
\end{enumerate}
We then compute, for each $\ell \in \N$ and $\sigma \in \Shadows$, the
Dirichlet generating function which enumerates the classes with shadow
$\sigma$, viz.\
$$ \zeta_\ell^\sigma(s) := \sum_{\mathcal{C} \in
\mathcal{Q}_\ell[\sigma]} \lvert \mathcal{C} \rvert^{-s}, \quad
\text{where $\mathcal{Q}_\ell[\sigma] := \{ \mathcal{C} \in
\GL_3(\lri)\backslash\mathfrak{gl}_3(\lri/\mathfrak{p}^\ell) \mid
\sh(\mathcal{C})=\sigma \}$.}
$$
Finally, we prove that if $p>3e(\lri\vert\Z_p)$ then
  \begin{equation}\label{zeta.shadows}
  \zeta_{\SL_3(\lri)}(s) = \lim_{\ell\to\infty} \lvert \lri/\mathfrak{p}
  \rvert^{-\ell} \sum_{\sigma \in \Shadows} \frac{ \lvert
    \GL_3(\lri/\mathfrak{p}) : H(\sigma) \rvert^{1-s/2} \,
    \zeta_{H(\sigma) \cap \SL_3(\lri/\mathfrak{p})}(s)} {\lvert
    \SL_3(\lri/\mathfrak{p}) : H(\sigma) \cap \SL_3(\lri/\mathfrak{p})
    \rvert} \, \zeta_\ell^{\sigma}(s/2),
  \end{equation}
  where, for each $\sigma\in\Shadows$, $H(\sigma) \leq
  \GL_3(\lri/\mathfrak{p})$ denotes a fixed representative of the
  conjugacy class $\sigma$. Theorem~\ref{zeta.SL3} is a consequence of
  \eqref{zeta.shadows}.  Theorem~\ref{thmABC:globalSL3} follows from
  an analysis of the formula given in Theorem~\ref{zeta.SL3}, and
  standard Tauberian theorems.

\begin{acknowledgements}
  The authors, in various constellations, would like to thank Alex
  Lubotzky and the following institutions for their support: the
  Batsheva de Rothschild Fund for the Advancement of Science, the
  EPSRC, the Mathematisches Forschungsinstitut Oberwolfach and the
  Nuffield Foundation.
\end{acknowledgements}

\end{document}